\documentclass{conm-p-l}

\usepackage[matrix,arrow,curve]{xy}

\usepackage[latin1]{inputenc}

\usepackage{amssymb}
\usepackage[mathscr]{eucal}

\newenvironment{env}[2]{\begin{#1}#2\end{#1}}{}
    \newcommand{\beq}[1]{\begin{env}{equation}{#1}}
    \newcommand{\beqn}[1]{\begin{env}{equation*}{#1}}
    \newcommand{\bal}[1]{\begin{env}{align}{#1}}
    \newcommand{\baln}[1]{\begin{env}{align*}{#1}}
    \newcommand{\bga}[1]{\begin{env}{gather}{#1}}
    \newcommand{\bgan}[1]{\begin{env}{gather*}{#1}}
    \newcommand{\bflal}[1]{\begin{env}{flalign}{#1}}
    \newcommand{\bflaln}[1]{\begin{env}{flalign*}{#1}}
    \newcommand{\bmu}[1]{\begin{env}{multline}{#1}}
    \newcommand{\bmun}[1]{\begin{env}{multline*}{#1}}
    \newcommand{\bsp}[1]{\begin{env}{split}{#1}}

    \newcommand{\eeq}{\end{env}}
    \newcommand{\eeqn}{\end{env}}
    \newcommand{\eal}{\end{env}}
    \newcommand{\ealn}{\end{env}}
    \newcommand{\ega}{\end{env}}
    \newcommand{\egan}{\end{env}}
    \newcommand{\eflal}{\end{env}}
    \newcommand{\eflaln}{\end{env}}
    \newcommand{\emu}{\end{env}}
    \newcommand{\emun}{\end{env}}
    \newcommand{\esp}{\end{env}}

\newcommand{\lf}{\vspace{2ex}}

\renewcommand{\bf}[1]{\textbf{#1}}
\renewcommand{\it}[1]{\textit{#1}}

\renewcommand{\sf}[1]{\textsf{#1}}

\newcommand{\hl}[1]{\bf{\it{#1}}}

\newcommand{\mbf}[1]{\mathbf{#1}}

\newcommand{\cmc}[1]{\mathcal{#1}}
\newcommand{\eus}[1]{\mathscr{#1}}
\newcommand{\euf}[1]{\mathfrak{#1}}
\newcommand{\bb}[1]{\mathbb{#1}}
\newcommand{\nbd}[1]{$#1$\nobreakdash--}
\newcommand{\ol}[1]{\overline{#1}}

\newcommand{\vt}{\vartheta}

\newcommand{\om}{\omega}
\newcommand{\Om}{\Omega}

\newcommand{\bfam}[1]{\bigl(#1\bigr)}
\newcommand{\Bfam}[1]{\Bigl(#1\Bigr)}
\newcommand{\AB}[1]{\langle#1\rangle}

\newcommand{\CB}[1]{\{#1\}}
\newcommand{\bCB}[1]{\bigl\{#1\bigr\}}
\newcommand{\BCB}[1]{\Bigl\{#1\Bigr\}}

\newcommand{\Matrix}[1]{\begin{pmatrix}#1\end{pmatrix}}

\newcommand{\sbars}[1]{\:\bar{#1}^s\:}
\newcommand{\set}[2][]{
    \ifthenelse{\equal{#1}{}}{
        \CB{#2}}{
        \CB{#1~|~#2}}}
\newcommand{\bset}[2][]{
    \ifthenelse{\equal{#1}{}}{
        \bCB{#2}}{
        \bCB{#1~|~#2}}}
\newcommand{\Bset}[2][]{
    \ifthenelse{\equal{#1}{}}{
        \BCB{#2}}{
        \BCB{#1~\big|~#2}}}

\DeclareMathOperator{\ls}{\normalfont\sf{span}}
\DeclareMathOperator{\cls}{\ol{\ls}}

\DeclareMathOperator{\id}{\normalfont\sf{id}}

\newcommand{\C}{\bb{C}}

\newcommand{\N}{\bb{N}}

\newcommand{\R}{\bb{R}}

\newcommand{\cA}{\cmc{A}}
\newcommand{\cB}{\cmc{B}}
\newcommand{\cC}{\cmc{C}}

\newcommand{\sB}{\eus{B}}

\newcommand{\sE}{\eus{E}}

\newcommand{\sK}{\eus{K}}

\newcommand{\ei}{\euf{i}}

\newcommand{\eH}{\euf{H}}

\newcommand{\U}{\mbf{1}}

\newcommand{\G}{\Gamma}

    \numberwithin{equation}{section}

        \newcommand{\definame}{Definition.}
        \newcommand{\propname}{Proposition.}
        \newcommand{\lemname}{Lemma.}
        \newcommand{\exname}{Example.}
        \newcommand{\remname}{Remark.}
        \newcommand{\obname}{Observation.}
        \newcommand{\corname}{Corollary.}

\swapnumbers
            \newtheorem{emp}{}[section]
                \newcommand{\bemp}[1][]{
                    \begin{emp}\hskip-\labelsep\bf{#1}\hskip\labelsep}
                \newcommand{\eemp}{\end{emp}}
            \newtheorem{itemp}[emp]{}
                \newcommand{\bitemp}[1][]{
                    \begin{itemp}\hskip-\labelsep\bf{#1}\hskip\labelsep\normalfont\itshape}
                \newcommand{\eitemp}{\end{itemp}}
            \newtheorem{thm}[emp]{\thmname}
                \newcommand{\bthm}{\begin{thm}}
                \newcommand{\ethm}{\end{thm}}
            \newtheorem{prop}[emp]{\propname}
                \newcommand{\bprop}{\begin{prop}}
                \newcommand{\eprop}{\end{prop}}
            \newtheorem{cor}[emp]{\corname}
                \newcommand{\bcor}{\begin{cor}}
                \newcommand{\ecor}{\end{cor}}
            \newtheorem{lem}[emp]{\lemname}
                \newcommand{\blem}{\begin{lem}}
                \newcommand{\elem}{\end{lem}}
    \theoremstyle{definition}
            \newtheorem{ex}[emp]{\exname}
                \newcommand{\bex}{\begin{ex}}
                \newcommand{\eex}{\end{ex}}
            \newtheorem{defi}[emp]{\definame}
                \newcommand{\bdefi}{\begin{defi}}
                \newcommand{\edefi}{\end{defi}}
    \theoremstyle{remark}
            \newtheorem{rem}[emp]{\remname}
                \newcommand{\brem}{\begin{rem}}
                \newcommand{\erem}{\end{rem}}
            \newtheorem{ob}[emp]{\obname}
                \newcommand{\bob}{\begin{ob}}
                \newcommand{\eob}{\end{ob}}

\begin{document}

\title[Topics Related to Product Systems]{Commutants of von Neumann Modules, Representations of $\sB^a(E)$ and Other Topics Related to Product Systems of Hilbert Modules}
\author{Michael Skeide}
\address{
Lehrstuhl f\"ur Wahrscheinlichkeitstheorie und Statistik\\
Brandenburgische Technische Universit\"at Cottbus\\
Postfach 10 13 44, D--03013 Cottbus, Germany
}
\curraddr{
Dipartimento S.E.G.e S.\\
Università degli Studi del Molise\\
Via de Sanctis\\
86100 Campobasso, Italy
}
\email{skeide@math.tu-cottbus.de}
\urladdr{www.math.tu-cottbus.de/INSTITUT/lswas/\_skeide.html}
\subjclass[2000]{Primary 46L55 60J25; Secondary 46L08 46L57 81S25}
\thanks{This work is supported by a PPP-project by DAAD and DST}





\begin{abstract}
We review some of our results from the theory of product systems of Hilbert modules \cite{BhSk00,BBLS00p,Ske00b,Ske01p,Ske02,Ske03b}. We explain that the product systems obtained from a CP-semigroup in \cite{BhSk00} and in \cite{MuSo02p} are commutants of each other. Then we use this new commutant technique to construct product systems from \nbd{E_0}semigroups on $\sB^a(E)$ where $E$ is a strongly full von Neumann module. (This improves the construction from \cite{Ske02} for Hilbert modules where existence of a unit vector is required.) Finally, we point out that the Arveson system of a CP-semigroup constructed by Powers from two spatial \nbd{E_0}semigroups is the product of the corresponding spatial Arveson systems as defined (for Hilbert modules) in \cite{Ske01p}. It need not coincide with the tensor product of Arveson systems.
\end{abstract}


\maketitle

$\sB(G)$ ($G$ some Hilbert space) are in many respects the simplest von Neumann algebras. They are factors which contain their commutant, they contain all finite rank operators and this determines completely the theory of all normal representations.

Hilbert modules (or, more precisely, von Neumann modules; see Section \ref{comsec}) $E$ over $\sB(G)$ share this simplicity. $E$ is always isomorphic to $\sB(G,H)$ ($H$ some other Hilbert space) equipped with the natural right module action and inner product $\AB{x,y}=x^*y$. The algebra $\sB^a(E)$ of all adjointable mappings on $E$ is just $\sB(H)$.

If $E$ is a two-sided \nbd{\sB(G)}module, i.e.\ if $H$ carries a normal representation of $\sB(G)$, then $H=G\otimes\eH$ ($\eH$ some other Hilbert space) and elements $b\in\sB(G)$ act in the natural way as $b\otimes\id_\eH$. The Hilbert space $\eH$ can be identified naturally with the space $\CB{x\in E\colon bx=xb~(b\in\sB(G))}$ of all elements in $E$ which intertwine left and right action of $\sB(G)$, where $h\in\eH$ corresponds to the intertwiner $\id_G\otimes h\colon g\mapsto g\otimes h$. (It is not difficult to see that every intertwiner arises in that way, and that $\AB{x,y}\in\sB(G)'=\C\U$ gives the correct scalar product as multiple of $\U$.) We observe that every von Neumann \nbd{\sB(G)}\nbd{\sB(G)}module $E=\sB(G)\sbars{\otimes}\eH$ ($\sbars{\otimes}$ stands for strong closure in $\sB(G,G\otimes\eH)$) contains its commutant with respect to $\sB(G)$ and is generated by it (in the strong topology) as a right module.

Let $\vt=\bfam{\vt_t}_{t\in\R_+}$ be a \hl{normal \nbd{E_0}semigroup} (a semigroup of normal unital endomorphisms) on $\sB(G)$. Arveson introduced product systems of Hilbert spaces (\hl{Arveson systems}, for short) in \cite{Arv89} by applying the preceding construction to the von Neumann $\sB(G)$--$\sB(G)$--mod\-ules $E_t=\sB(G)$ with left action $b.x_t=\vt_t(b)x_t$. Thus, $E=\sB(G,G\otimes\eH_t)$ and it is not difficult to see that
\beqn{
\eH_{s+t}
~=~
\eH_s\otimes\eH_t
}\eeqn
in a natural (i.e.\ associative) way. Factorizable families of Hilbert spaces have been recognized before (Streater \cite{Str69}, Araki \cite{Ara70}, Parthasarathy and Schmidt \cite{PaSchm72}). However, the examples were always of the simplest possible type I (symmetric Fock spaces). The systematic study of Arveson systems started with \cite{Arv89}.

A different approach to Arveson systems was discovered by Bhat \cite{Bha96}, by inventing an elegant and quick proof for the representation theory of $\sB(G)$. Indeed, if we choose a unit vector $\Om\in G$, if we denote by $p_t=\vt_t(\Om\Om^*)$ the time evolution of the rank-one projection $\Om\Om^*$, and if we set $G_t=p_tG$, then it is not difficult to see that
\beqn{
g\otimes g_t
~\longmapsto~
\vt_t(g\Om^*)g_t
}\eeqn
defines a unitary isomorphism $G\otimes G_t\rightarrow G$. (We repeat this more detailed for Hilbert modules in Section \ref{pssec}.) The restriction of this isomorphism gives isomorphisms $G_s\otimes G_t=G_{s+t}$ and it is not difficult to see that the Arveson systems $\eH^\otimes=\bfam{\eH_t}_{t\in\R_+}$ and $G^\otimes=\bfam{G_t}_{t\in\R_+}$ are (anti-) isomorphic.

The construction of Arveson systems from \nbd{E_0}semigroups has been generalized by Bhat \cite{Bha96} and Arveson \cite{Arv97} to \hl{normal CP-semigroups} (semigroups of (unital) normal completely positive mappings) on $\sB(G)$ with the help of minimal weak dilations. A \hl{weak dilation} of a normal CP-semigroup $T=\bfam{T_t}_{t\in\R_+}$ on $\sB(G)$ is a normal \nbd{E_0}semigroup $\vt=\bfam{\vt_t}_{t\in\R_+}$ on $\sB(H)$, where $H\supset G$ is another Hilbert space, such that the diagram
\beqn{
\parbox{8cm}{
\xymatrix{
\sB(G)	\ar[rr]^{T_t}	\ar[d]_{\ei}	&&	\sB(G)
\\
\sB(H)	\ar[rr]_{\vt_t}			&&	\sB(H)	\ar[u]_{p\bullet p}
}
}
}\eeqn
commutes, where $\ei$ is the natural embedding of $\sB(G)$ as a corner into $\sB(H)$, and where $p$ is the projection onto $G$. The dilation is \hl{minimal}, if the vectors $\vt_t\circ\ei(b)g$ ($t\in\R_+,b\in\sB(G),g\in G$) are total in $H$. The minimal weak dilation is a unique universal object in the category of weak dilations and, therefore, there is no arbitrariness in saying the Arveson system of $T$ is the Arveson system of the unique minimal dilating \nbd{E_0}semigroup $\vt$.

The intersection of interest of quantum dynamical systems and quantum probability is dilation theory. CP-semigroups on $\sB(G)$ can be dilated in many ways to \nbd{E_0}semigroups on $\sB(H)$. Therefore, as long as we are interested in systems described by a full algebra $\sB(G)$, dilations can be delt with on $\sB(H)$. As soon as we want to dilate CP-semigroups on a more general unital \nbd{C^*}algebra $\cB$ (e.g.\ in order to include also classical dynamical systems, where all algebras are commutative), this is no longer true. For instance, Bhat \cite{Bha99} constructed for every such semigroup a dilation to a non-unital \nbd{C^*}subalgebra $\cA$ of some $\sB(H)$ and it was not possible to extend the dilating \nbd{E_0}semigroup to all of $\sB(H)$.

It turns out that a good intermediate object on which dilating \nbd{E_0}semigroups can act is the algebra $\sB^a(E)$ of adjointable mappings on a Hilbert module over the algebra $\cB$. As explicit constructions of dilations from Bhat and Skeide \cite{BhSk00} and Skeide \cite{Ske00} show, $\sB^a(E)$ is small enough ($\vt$ as constructed in \cite{Bha99} extends from $\cA\subset\sB^a(E)\subset\sB(H)$ to $\sB^a(E)\subset\sB(H)$, but not to $\sB(H)$), but, still sufficiently simple ($\sB^a(E)$ is generated by its rank-one operators). We have the following situation
\beq{\label{dildig}\tag{$*$}
\parbox{8cm}{
\xymatrix{
\cB	\ar[rr]^{T_t}	\ar[d]_{\xi\bullet\xi^*}	&&	\cB
\\
\sB^a(E)	\ar[rr]_{\vt_t}			&&	\sB^a(E)	\ar[u]_{\AB{\xi,\bullet\xi}}
}
}
}\eeq
of weak dilation on a Hilbert module, where $\xi$ is a unit vector (i.e.\ $\AB{\xi,\xi}=\U$ and we say $E$ is unital).

It is the goal of these notes to explain some of the possibilites. In Section \ref{pssec} we describe our construction from \cite{Ske02} of product systems of Hilbert modules (\hl{product systems}, for short) from \nbd{E_0}semigroups on $\sB^a(E)$, in the case when $E$ is unital. This construction is the analogue of Bhat's \cite{Bha96} approach to Arveson systems. It provides a complete treatment of the theory of strict representations of $\sB^a(E)$ on another Hilbert module for unital Hilbert modules $E$.

The unit vector plays the same role in both cases and the we use it underlines the importance of the fact that both $\sB(H)$ and $\sB^a(E)$ are generated by their rank-one operaotrs. If $E$ is the module on which a dilation acts, then the requirment for existence of a unit vector (a non-trivial requirment, in general; see Example \ref{non1ex}) is always fulfilled. Nevertheless, it remains the question, whether it is possible to do the representation theory also for Hilbert modules without a unit vector. In Section \ref{comsec} we show (for the first time) that this is possible for normal representations of $\sB^a(E)$ when $E$ is a full von Neumann module. The proof is completely different from that in \cite{Ske02}. It is inspired very much by an idea from Muhly and Solel \cite{MuSo02p} (see also their contribution to this meeting) to construct two-sided von Neumann modules over the commutant of $\cB$. The main new ingredient is to establish a one-to-one correspondence between von Neumann \nbd{\cB}\nbd{\cB}modules and von Neumann \nbd{\cB'}\nbd{\cB'}modules which generalizes the relation between $\cB$ and $\cB'$.

Product systems of Hilbert modules have appeared first probably in A.\ Alevras' PhD-thesis (reference unkown) in an attempt to construct product systems from \nbd{E_0}semigroups on type II factors in analogy with Arveson's intertwiner spaces. See also the contributions of R.\ Floricel and I.\ Hirshberg to this meeting. The first product systems of two-sided Hilbert modules over general \nbd{C^*}algebras appeared probably in Bhat and Skeide \cite{BhSk00} where we constructed a product system directly from a CP-semigroup. In the case of CP-semigroups on $\sB(G)$ we are able to construct directly the associated Arveson system without the way arround via minimal weak dilation. In Section \ref{powsec} we use this possibility to answer a question raised by R.\ Powers in this meeting, whether a certain CP-semigroup constructed from two spatial \nbd{E_0}semigroups must have an Arveson system which is the tensor product of the Arveson systems of the \nbd{E_0}semigroups in the negative sense.

Finally, we explain that the product system of von Neumann \nbd{\cB'}\nbd{\cB'}modules constructed in \cite{MuSo02p} from a CP-semigroup on a von Neumann algebra $\cB$ is just the commutant of the product system of von Neumann \nbd{\cB}\nbd{\cB}modules as constructed on \cite{BhSk00}. We are able to generalize their construction to a construction starting from an \nbd{E_0}semigroup on $\sB^a(E)$. This construction of product systems from \nbd{E_0}semigroups on $\sB^a(E)$ (based on intertwiner spaces) yields the commutant system of that obtained in our generalization \cite{Ske02} of Bhat's approach to Arveson systems (based on the representation theory of $\sB(H)$).

\section{Product systems and dilations of CP-semigroups} \label{pssec}

In this section we review some of our results about product systems in dilation theory. Some results about classification of product systems can be found in the beginning of Section \ref{powsec}. The reader who wishes a more complete account about product systems might consult the survey Skeide \cite{Ske03b}.

Let $\cB$ be a unital \nbd{C^*}algebra and let $E$ be a Hilbert \nbd{\cB}module. The \hl{strict topology} on $\sB^a(E)$, the algebra of all mappings $a\colon E\rightarrow E$ which have an adjoint, is that obtained by considering $\sB^a(E)$ as mutliplier algebra of $\sK(E)$, the \nbd{C^*}algebra of \hl{compact operators} (i.e.\ the norm closure of the \nbd{*}algebra generated by the \hl{rank-one operators} $xy^*\colon z\mapsto x\AB{y,z}$). On bounded subsets the strict topology coincides with the \nbd{*}strong topology of $\sB^a(E)$.

Suppose $E$ is \hl{unital}, i.e.\ $E$ has a \hl{unit vector} $\xi$, and let $\vt=\bfam{\vt_t}_{t\in\R_+}$ be a strict \nbd{E_0}semigroup on $\sB^a(E)$, i.e.\ all mappings $\vt_t$ are continuous in the strict topology. Denote by $p_t=\vt_t(\xi\xi^*)$ the time evolution of the rank-one projection $\xi\xi^*$. Then $E_t=p_tE$ is a Hilbert \nbd{\cB}submodule of $E$. It is easy to check that $\ei\colon b\mapsto\xi b\xi^*$ defines a representation of $\cB$ on $E$ (cf.\ Diagram \eqref{dildig}). We turn $E_t$ into a two-sided Hilbert \nbd{\cB}module by defining the left action $bx_t=\vt_t\circ\ei(b)x_t$. Observe that this left action is unital. Denoting by $\odot$ the \hl{tensor product over $\cB$}, the computation
\beqn{
\AB{x\odot y_t,x'\odot y'_t}
~=~
\AB{y_t,\AB{x,x'}y_t}
~=~
\AB{y_t,\vt_t(\xi x^*x'\xi^*)y_t}
~=~
\AB{\vt_t(x\xi^*)y_t,\vt_t(x'\xi^*)y_t}
}\eeqn
shows that $x\odot y_t\mapsto \vt_t(x\xi^*)y_t$ defines an isometry $u_t\colon E\odot E_t\rightarrow E$. Moreover, since $\sK(E)$ has an approximate unit $\Bfam{\sum\limits_{k=1}^{n_\lambda}v_k^\lambda{w_k^\lambda}^*}_\lambda$, and since $\vt_t$ is strict, we find that
\beqn{
x
~=~
\lim_\lambda\vt_t\Bfam{\sum\limits_{k=1}^{n_\lambda}v_k^\lambda{w_k^\lambda}^*}x
~=~
\lim_\lambda\sum\limits_{k=1}^{n_\lambda}\vt_t(v_k^\lambda\xi^*)\vt_t(\xi{w_k^\lambda}^*)x
~\in~
u_t(E\odot E_t),
}\eeqn
(because $\vt_t(\xi{w_k^\lambda}^*)x=p_t\vt_t(\xi{w_k^\lambda}^*)x\in E_t$,) so that $u_t$ is a unitary. We recover $\vt_t$ as $\vt_t(a)=a\odot\id_{E_t}$. (Of course, this concerns so far only a single endomorphism $\vt_t$, and, actually, it is not important that $\vt_t$ maps into $\sB^a(E)$. We find the same results for aribitrary strict representations of $\sB^a(E)$ on another Hilbert \nbd{\cC}module $F$; see \cite{Ske02}. One may, for instance, show that the unital Hilbert modules $E$ and $F$ have strictly isomorphic operator algebras, if and only if there is a Morita equivalence \nbd{\cB}\nbd{\cC}module $F_\xi$ such that $F=E\odot F_\xi$; see \cite{MSS02p}.)

It is not difficult to see that the restriction $u_{st}$ of $u_t$ to $E_s\odot E_t$ defines a \hl{two-sided} (i.e.\ a \nbd{\cB}\nbd{\cB}linear) unitary $E_s\odot E_t\rightarrow E_{s+t}$. Therefore,
\beqn{
E_s\odot E_t
~=~
E_{s+t}
}\eeqn
(associatively) and $E^\odot=\bfam{E_t}_{t\in\R_+}$ is a product system in the sense of \cite{BhSk00}. One of the deepest results in Arveson's theory asserts that every Arveson system comes from an \nbd{E_0}semigroup on $\sB(H)$. Observe that the analogue problem for Hilbert modules is completely open (except the case, when there exist units, which is easy; see below, when we discuss dilations). Like in Arveson's theory \nbd{E_0}semigroups on the same $\sB^a(E)$ are determined by their product system up to cocycle conjugacy. We do not yet know, whether there exists an analogue for \nbd{E_0}semigroups on $\sB^a(E)$ and $\sB^a(E')$ when $E$ and $E'$ are not necessarily isomorphic. (In Arveson's theory $\sB(H)\cong\sB(H')$ is a hidden assumption, because all Hilbert spaces are infinite-dimensional and separable.)

We ask under which circumstances the \nbd{E_0}semigroup $\vt$ and the unit vector $\xi$ define a \hl{weak dilation} (see Diagram \eqref{dildig}). In other words, when does $T_t(b)=\AB{\xi,\vt_t\circ\ei(b)\xi}$ define a CP-semigroup on $\cB$? It is easy to check that this happens, if and only if $p_t\ge p_0$ for all $t$. In this case, we have $\xi_t:=\xi\in E_t$ for all $t$ and
\beqn{
\xi_s\odot\xi_t
~=~
\xi_{s+t}
}\eeqn
so that $\xi^\odot=\bfam{\xi_t}_{t\in\R_+}$ is a \hl{unit} for $E^\odot$ in the sense of \cite{BhSk00} with $T_t(b)=\AB{\xi_t,b\xi_t}$.

Conversely, if $E^\odot$ is a product system with a (unital) unit $\xi^\odot$, then $T_t(b)=\AB{\xi_t,b\xi_t}$ defines a (unital) CP-semigroup, and by an inductive limit construction in \cite{BhSk00} (embedding $E_t$ as $\xi_s\odot E_t$ into $E_{s+t}$) we construct a Hilbert module $E'$ with an \nbd{E_0}semigroup dilating $T$. If in the preceding paragraph the projections $p_t$ increase to $\U$, then $E'$ is the module $E$ we started with.

Finally, by another inductive limit construction in \cite{BhSk00} we can construct from every CP-semigroup $T$ a (unique) product system $E^\odot$ and a unit $\xi^\odot$, such that $E^\odot$ is generated by $\xi^\odot$ (there is no proper subsystem containing $\xi^\odot$), and such that $T_t(b)=\AB{\xi_t,b\xi_t}$.

For a more detailed introduction to these last constructions we refer the reader to the survey \cite{Ske03b} and, of course, to the original paper \cite{BhSk00}.

\section{Von Neumann modules, commutants and endomorphisms of $\sB^a(E)$} \label{comsec}

In this section we repeat the definition of von Neumann modules  \cite{Ske00b}. Starting from the fact that a von Neumann algebra $\cB$ is a concrete subalgebra of operators on a Hilbert space $G$, the definition of von Neumann modules makes use of the possibility to construct an explicit identification of a Hilbert \nbd{\cB}module as a subspace of $\sB(G,H)$ where $H$ is another canonically associated Hilbert space. Then we show, based on a result from Muhly and Solel \cite{MuSo02p}, how it is possible to associate with a two-sided von Neumann module \nbd{\cB}module $E$ (a \nbd{\cB}correspondence in the terminology of \cite{MuSo02p}) a \hl{commutant} $E'$ which is a two-sided von Neumann module over the commutant $\cB'$ of $\cB$ in $\sB(G)$. Thanks to speaking about von Neumann modules, i.e.\ thinking always about concrete subspaces of $\sB(G)$ and $\sB(G,H)$, all the occuring functors are one-to-one (and not only up to isomorphism as in \cite{MuSo02p}). We explain that the product systems constructed in \cite{MuSo02p} and in in \cite{BhSk00} are just commutants of each other. Finally, we sketch how to use the idea of commutant to determine completely the theory of normal representations of $\sB^a(E)$ on a von Neumann \nbd{\cB}module $F$ for an arbitrary von Neumann \nbd{\cC}module $E$. This improves the construction from Section \ref{pssec} where we required existence of a unit vector. We use the representation theory to construct a product system from a normal \nbd{E_0}semigroup on $\sB^a(E)$ also without existence of a unit vector in $E$.

Let $\cB\subset\sB(G)$ be a von Neumann algebra on a Hilbert space $G$. (Now and in the sequel, all representations, also the defining one, are assumed to act non-degenerately.)  Let $E$ be a Hilbert \nbd{\cB}module and set $H=E\odot G$. This is the module tensor product of a Hilbert \nbd{\cB}module and a Hilbert \nbd{\cB}\nbd{\C}module and, therefore, a Hilbert space. For $x\in E$ define $L_x\in\sB(G,H)$ by setting $L_xg=x\odot g$. Then
\baln{
\AB{x,y}
~=~
L_x^*L_y
&&&
L_{xb}
~=~
L_xb.
}\ealn
Indentifying $L_x=x$ we find $E\subset\sB(G,H)$. We say $E$ is a \hl{von Neumann \nbd{\cB}module}, if $E$ is strongly closed in $\sB(G,H)$. One may show that $E$ is a von Neumann module, if and only if $E$ is self-dual. Another important result is that submodules of $E$ with zero-complement are strongly dense in $E$.

$H$ carries the unital faithful representation $\pi(a)\colon x\odot g\mapsto ax\odot g$ of $\sB^a(E)$ and the range of $\pi$ is a von Neumann subalgebra of $\sB(H)$ which we identify with $\sB^a(E)$. If $\cA$ is another von Neumann algebra and $E$ also a Hilbert \nbd{\cA}\nbd{\cB}module, then $E$ is a \hl{von Neumann \nbd{\cA}\nbd{\cB}module}, if the representation $\cA\rightarrow\sB^a(E)\rightarrow\sB(H)$ is normal.

\brem
If $\xi\in E$ such that $\cls^s\cA\xi\cB$, then $(E,\xi)$ is the GNS-con\-struc\-tion \cite{Pas73} for the CP-map $a\mapsto\AB{\xi,a\xi}$, while $a\mapsto L_\xi^*\pi(a)L_\xi$ is the Stinespring construction.
\erem

$\rho'(b')\colon x\odot g\mapsto x\odot b'g$ defines a normal unital representation $\rho'$ of $\cB'$ on $H$. Straightforward computations show that the intertwiner space
\beqn{
C_{\cB'}(\sB(G,H))
~=~
\bCB{x\in\sB(G,H)\colon\rho'(b')x=xb' ~(b'\in\cB')}
}\eeqn
is a von Neumann \nbd{\cB}module containing $E$. (It is a right \nbd{\cB}module and the inner product $\AB{x,y}=x^*y$ has values in $\cB$.) Actually, $C_{\cB'}(\sB(G,H))=E$, because $\cls EG=H$ (i.e.\ $E$ has zero-complement and, therefore, is strongly dense). The following result from Muhly and Solel \cite{MuSo02p} shows, in a sense, the converse.

\bitemp[Lemma \cite{MuSo02p}.~]\label{MuSolem}
If $\rho'$ is a normal representation of $\cB'$ on a Hilbert space $H$, then
\beqn{
\cls C_{\cB'}(\sB(G,H))G
~=~
H.
}\eeqn
\eitemp

\noindent
Thus, we have a 1--1--correspondence between von Neumann \nbd{\cB}modules and representations of $\cB'$ (established by sending $(\rho',H)$ to $C_{\cB'}(\sB(G,H))$ and $E$ to $(\rho',E\odot G)$).

If $E$ is a von Neumann \nbd{\cB}\nbd{\cB}module, then $\rho=\pi\upharpoonright\cB$ is a representation of $\cB$ ``commuting'' with $\rho'$, i.e.\ $\rho(b)\rho'(b')=\rho'(b')\rho(b)$ ($b\in\cB,b'\in\cB'$). As the tripel $(H,\rho,\rho')$ is symmetric in $\cB\leftrightarrow\cB'$ we find equivalences
\beqn{
E
~\longleftrightarrow~
(H,\rho,\rho')
~\longleftrightarrow~
E'
}\eeqn
($E':=C_\cB(\sB(G,H))$) between von Neumann \nbd{\cB}\nbd{\cB}modules, triples $(H,\rho,\rho')$ and von Neumann \nbd{\cB'}\nbd{\cB'}modules (with bilinear mappings as morphisms).

\bitemp[Theorem (more or less \cite{MuSo02p}).]
\hfill
$(E_1\odot E_2)'~=~E'_2\odot E'_1$.
\hfill
~~~~~~
\eitemp

\begin{proof}
(Sketch.) By multiple use of Lemma \ref{MuSolem} and associativity of the tensor product we have
\beqn{
E_i\odot G=E'_i\odot G~~\text{and}~~E'_2\odot E'_1\odot G=E'_2\odot(E_1\odot G)=E_1\odot(E'_2\odot G)=E_1\odot E_2\odot G,
}\eeqn
where the only tricky identification $E'_2\odot(E_1\odot G)=E_1\odot(E'_2\odot G)$ is, indeed, done by $x'_2\odot y_1\odot g\mapsto y_1\odot x'_2\odot g$. The remaining things follow by looking at how the various representations of $\cB$ and $\cB'$ are defined.
\end{proof}

\lf\noindent
So product systems of \nbd{\cB}\nbd{\cB}modules (anti-)correspond to such of \nbd{\cB'}\nbd{\cB'}modules. In particular, product systems of \nbd{\sB(G)}\nbd{\sB(G)}modules (anti-)correspond to Arveson systems. It is easy to see directly from the definition of the product system in \cite{MuSo02p} that it is the commutant of the product system constructed in \cite{BhSk00}.

A von Neumann \nbd{\cB}module $E$ is \hl{strongly full}, if $\cB_E:=\cls^s\AB{E,E}=\cB$. We may always pass to a full von Neumann module restricting $\cB$ to $\cB_E$. However, in general, we may not assume that $E$ has a unit vector.

\bex\label{non1ex}
Let $\cB=\Bfam{\substack{\C\\[1ex]0}~\substack{0\\[1ex]\smash{M_2}}}\subset M_3$. Then $E=\Bfam{\!\!\substack{~~0~~{\C^2}^*\\[.5ex]{\C^2}~0~}\!}\subset M_3$ is a Morita equivalence $\cB$--$\cB$--module without a unit vector. In particular, $E$ is a full non-unital von Neumann module.
\eex

We present now a new construction of a product system from a normal $E_0$--semigroup $\vt$ on $\sB^a(E)$ for an arbitrary full von Neumann \nbd{\cB}module $E$. It is based on the observation that the commutant of $\sB^a(E)$ in $\sB(H)$ is just $\rho'(\cB')$. (This can be shown most easily --- although not most elementarilly --- by considering $E$ as Morita equivalence module between the von Neumann algebras $\cB$ and $\sB^a(E)$.) $\vt_t$ induces a second representation $\rho_t=\pi\circ\vt_t$ of $\sB^a(E)$ on $H=E\odot G$. Set $E'_t=\bCB{x'_t\in\sB(H)\colon\rho_t(a)x'_t=x'_ta~(a\in\sB^a(E))}$. Then $E'_t$ is a von Neumann bimodule over $\sB^a(E)'=\rho'(\cB')$. Since $E$ is full, $\rho'$ is faithful. So we may interprete $E'_t$ as von Neumann \nbd{\cB'}\nbd{\cB'}module. From the identification $E'_t\subset\sB(H)$ it follows that $E'_s\odot E'_t$ embeds into $E'_{s+t}$ and, once again, from the crucial Lemma \ref{MuSolem} it follows that this embedding is surjective. (Associativity is gifted by that of multiplication in $\sB(H)$.) Therefore, ${E'}^\odot=\bfam{E'_t}_{t\in\R_+}$ is a product system of \nbd{\cB'}\nbd{\cB'}modules and, consequently, $E^\odot=\bfam{E_t}_{t\in\R_+}$, defined by setting $E_t=(E'_t)'$, is a product system of \nbd{\cB}\nbd{\cB}modules.

Most important are the following equalities
\beqn{
E
~=~
E\odot E_t
~~~~~~~~~
\vt_t(a)
~=~
a\odot\id_{E_t}
}\eeqn
which settle (as the corresponding formulae in Section \ref{pssec}) the representation theory of $\sB^a(E)$. The identification is done by sending $x\odot y'_t\odot g\in E\odot(E'_t\odot G)=E\odot E_t\odot G$ to $y'_t(x\odot g)\in E\odot G$. Clearly, this identification intertwines the corresponding representations of $\cB'$ so that, indeed, the modules (being intertwiners for the same representation) coincide. Since $y'_t$, by definition, intertwines $\vt_t(a)$ and $a$ also the second equality follows. The preceding constructions are joint work with P.\ Muhly and B.\ Solel. Detailed proofs and some more results will appear in \cite{MSS02p}. For instance, we will show that the question whether $E^\odot$ comes from an \nbd{E_0}semigroup $\vt$ is equivalent to the question whether ${E'}^\odot$ admits a unitary representation on some $\sB(H)$ in the sense of \cite{MuSo02p}.

\section{The product of spatial product systems and a problem of Bob Powers} \label{powsec}

The classification of Arveson systems, in a first step, is based on the existence of units. If the Arveson system is generated by its units, then it is type I, if there exist units, but they do not generate the system, then it is type II, and if there are no units, then it is type III. Non-type III Arveson systems are also called spatial. Spatial Arveson systems contain a maximal type I (or completely spatial) subsystem. Type I Arveson systems are symmetric Fock spaces and as such characterized by a single Hilbert space, i.e.\ more or less by a dimension. The dimension of this space for the maximal completely spatial subsystem of a spatial Arveson systems is its Arveson index. It is possible to build the tensor product of Arveson systems and the index behaves additive under tensor product.

For Hilbert modules we have many parallels, but, sometimes the situation is more complicated. Also here we make a first distinction into type I, type II and type III accordings to how many units exist. (We should say that Arveson units are measurable, hence, ``continuous'' in a suitable sense. We speak only about ``continuous'' units; see \cite{Ske03b} for details.) The analogue of the symmetric Fock space is the time-ordered Fock module; see \cite{BhSk00,LiSk01}. Unfortunately, it is easily possible to write down examples for type I product systems which are not time-ordered Fock modules (see \cite{Ske03b}). On the other hand, type I systems of von Neumann modules are always time-ordered Fock modules \cite{BBLS00p}. The crucial point is to show existence of a \hl{central} unit, i.e.\ a unit whose elements commute with all algebra elements. (For Hilbert spaces this is trivial, because every vector commutes with $\C$. For von Neumann modules the result from \cite{BBLS00p} that a product system with units has also central units is, in fact, equivalent to the results by Christensen and Evans \cite{ChrEv79} on the generator of a normal uniformly continuous CP-semigroup on a von Neumann algebra.) Once a central unit is given, a type I system must be a Fock module even in the case of Hilbert modules. Also here a product system with a central unit contains a maximal type I subsystem which is a Fock module. In Skeide \cite{Ske01p} we called such product systems \hl{spatial} and defined the \hl{index} as the bimodule which characterizes the Fock module. (This bimodule is no longer determined by a simple dimension.) Summarizing, for spatial product systems of Hilbert modules we have analogy with Arveson systems, and product systems of von Neumann modules, if not type III, are spatial automatically.

Another thing which does not work, is the tensor product of product systems. (At least not, if we insist to stay inside the category of \nbd{\cB}modules. Of course, it is possible to consider the exterior tensor product of product systems.) In \cite{Ske01p} we constructed a product of spatial product systems, roughly speaking, by taking the central units as reference unit which are identified in the product. The product contains the original product systems and is generated by them. Components from different factors which are orthogonal to the respective reference units are orthogonal to each other in the product. By these properties the product $(E^1\circledcirc E^2)^\odot$ of spatial product systems ${E^i}^\odot$ with central (and unital) reference units ${\om^i}^\odot$ is already determined uniquely. More precisely, $(E^1\circledcirc E^2)_t$ is spanned by expressions
\beq{\label{pildef}
x_n\odot\ldots\odot x_1
~~~~~~~(n\in\N\,,\,t_n+\ldots+t_1=t\,,\,x_i\in E^1_{t_i}\vee x_i\in E^2_{t_i})
}\eeq
with inner product
\beq{\label{tipdef}
\AB{x_t,y_t}
~=~
\begin{cases}
\AB{x_t,y_t}&\text{for~}x_t,y_t\in E^1_t\vee x_t\,,\,y_t\in E^2_t
\\
\AB{x_t,\om^1_t}\AB{\om^2_t,y_t}&\text{for~}x_t\in E^1_t\,,\,y_t\in E^2_t
\\
\AB{x_t,\om^2_t}\AB{\om^1_t,y_t}&\text{for~}x_t\in E^2_t\,,\,y_t\in E^1_t.
\end{cases}
}\eeq
It can be shown by an inductive limit that such an object exists. The index (i.e.\ the module which determines the Fock module isomorphic to the maximal type I subsystem) behaves additively (direct sum) under the product. In the case of spatial Arveson systems our product is a subsystem of the tensor product. If one of the factors is type I, then it coincides with the tensor product. However, one can show (see Liebscher \cite{Lie00p1}) that there are type II Arveson systems whose product is different from their tensor product.

In this meeting R.\ Powers constructed a CP-semigroup $T$ on $\sB(G\oplus G)$ starting from two \nbd{E_0}semigroups $\vt^i$ on $\sB(G)$ with spatial product systems ${\eH^i}^\otimes$ and reference units ${\om^i}^\otimes$ by setting
\beqn{
T_t
\Matrix{
a_{11}&a_{12}
\\
a_{21}&a_{22}
}
~=~
\Matrix{
a_{11}\otimes\id_{\eH^1_t}&(\id_G\otimes\om^1_t)a_{12}(\id_G\otimes\om^2_t)^*
\\
(\id_G\otimes\om^2_t)a_{21}(\id_G\otimes\om^1_t)^*&a_{22}\otimes\id_{\eH^2_t}
}.
}\eeqn
(We make use of the identifications $G\otimes\eH^1_t=G=G\otimes\eH^2_t$. For instance, $a_{ii}\otimes\id_{\eH^i_t}$ is just $\vt^i_t(a_{ii})$.) Powers asked for the Arveson system of this CP-semigroup and, in particular, whether it could be the tensor product. In the remainder we will outline a proof that the Arveson system is just our product of the Arveson systems of $\vt^i$. Together with the preceding remark this shows that the Arveson system of $T$ can be the tensor product, but, there are counter examples.

One problem Powers had to face is that, in the usual approach, one first has to find the minimal dilation of $T$ and then to find the Arveson system of the dilating \nbd{E_0}semigroup. Our approach from \cite{BhSk00}, instead, is more direct. We construct directly from a CP-semigroup a product system $E^\odot$ of two-sided Hilbert modules (in this case two-sided von Neumann \nbd{\sB(G\oplus G)}modules). As outlined in the introduction, as two-sided von Neumann modules over $\sB(G\oplus G)$ the $E_t$ have the form $E_t=\sB(G\oplus G)\sbars{\otimes}\eH_t=\sB\bfam{\bfam{\substack{G\\G}},\bfam{\substack{G\\G}}\otimes\eH_t}$ and it is not difficult to understand (see \cite{BhSk00}) that $\eH^\otimes=\bfam{\eH_t}_{t\in\R_+}$ is nothing but the Arveson system of $T$.

The members $E_t$ of the product system from \cite{BhSk00} are constructed as inductive limits over
\beqn{
x^n_{t_n}\odot\ldots\odot x^1_{t_1}~~~~~~(n\in\N\,,\,t_n+\ldots+t_1=t\,,\,x^i_{t_i}\in\sE_{t_i})
}\eeqn
where $\sE_t$ is the GNS-module of $T_t$, i.e.\ $\sE_t=\cls^s\cB\xi_t\cB$ and $\AB{\xi_t,b\xi_t}=T_t(b)$. By comparison with \eqref{pildef} it is almost clear that we are done, if we find back in the GNS-construction for $T_t$ the inner product structure from \eqref{tipdef}.

We claim $\sE_t=\sB\bfam{\bfam{\substack{G\\G}},\bfam{\substack{G\\G}}\otimes H_t}$ with $H_t=\C\om_t\oplus(\eH^1_t\ominus\C\om^1_t)\oplus(\eH^2_t\ominus\C\om^2_t)$ and with the cyclic vector $\xi_t$ explained in the following way. Let $\bfam{\substack{g^1\\g^2}}=\bfam{\substack{h^1\otimes h^1_t\\h^2\otimes h^2_t}}\in\bfam{\substack{G\\G}}=\bfam{\substack{G\otimes\eH^1_t\\G\otimes\eH^2_t}}$. Denote $p^i_t=\id_{\eH^i_t}-\om^i_t{\om^i_t}^*$. Then put
\beqn{
\xi_t
\Matrix{
g^1\\g^2
}
~=~
\Matrix{h_1\\0}\otimes\bfam{\AB{\om^1_t,h^1_t}\,,\,p^1_th^1_t\,,\,0}
+
\Matrix{0\\h_2}\otimes\bfam{\AB{\om^2_t,h^2_t}\,,\,0\,,\,p^2_th^2_t}.
}\eeqn
One easliy veryfies $\AB{\xi_t,a\xi_t}=T_t(a)$ and also that $\xi_t$ generates $\sE_t$ as two-sided von Neumann module. Obviously, inner products of elements $h^i_t\in\eH^i_t\subset H_t$ and $h^j_t\in\eH^j_t\subset H_t$ behave exactly as required by \eqref{tipdef}. This shows that the Arveson system $\eH^\otimes$ is exactly $(\eH^1\circledcirc\eH^2)^\otimes$.

We mention that we are also able to write down the minimal dilation of $T$.

We close with the remark that the whole construction also works, if we replace $G\oplus G$ with $E^1\oplus E^2$ for Hilbert \nbd{\cB}modules $E^i$. Given two \nbd{E_0}semigroups $\vt^i$ on $\sB^a(E^i)$ with spatial product systems we can define $T$ in a similar manner, and the product system of $T$ is the product of the product systems of $\vt^i$. The preceding construction and the generalization to Hilbert modules are joint work with B.V.R.\ Bhat and V.\ Liebscher and will appear together with detailed proofs in \cite{BLS02p}.

\setlength{\baselineskip}{2.5ex}


\newcommand{\Swap}[2]{#2#1}\newcommand{\Sort}[1]{}
\providecommand{\bysame}{\leavevmode\hbox to3em{\hrulefill}\thinspace}
\providecommand{\MR}{\relax\ifhmode\unskip\space\fi MR }
\providecommand{\MRhref}[2]{%
  \href{http://www.ams.org/mathscinet-getitem?mr=#1}{#2}
}
\providecommand{\href}[2]{#2}

\end{document}